\documentclass[12pt]{article}       
\usepackage{a4}
\usepackage{amsmath,amssymb,amsthm} 
\usepackage[all]{xy}
\theoremstyle{plain}               
\newtheorem{defn}{Definition}[section]

\newtheorem{proposition}[defn]{Proposition}
\newtheorem{corollary}[defn]{Corollary}

\theoremstyle{definition}          
\newtheorem{remark}[defn]{Remark}
\newtheorem{example}[defn]{Example}
\begin{document}
\title{Planar Binomial Coefficients}
\author{Lothar Gerritzen}
\date{06.01.2005}
\maketitle

\begin{abstract}
The notion of binomial coefficients $T \choose S$ of finite planar, reduced rooted trees $T, S$ is defined and a
recursive formula for its computation is shown. The nonassociative binomial formula\\
$$(1 + x)^T = \displaystyle \sum_S  {T \choose S}  x^S$$\\
for powers relative to $T$ is derived.\\
Similarly binomial coefficients $ T \choose S, V$ of the second kind are introduced and it is shown that
$(x \otimes 1 + 1 \otimes x)^T= \displaystyle \sum_{S, V} {T \choose S, V} (x^S \otimes x^V)$\\
The roots $\sqrt[T]{1+x}= (1 + x) ^{T^{-1}}$ which are planar
power series $f$ such that $ f^T= 1+x$ are considered. Formulas
for their coefficients are given.
\end{abstract}

\section{Introduction}
The binomial coefficient ${T \choose S}$ of two finite, planar, reduced rooted trees $T, S$ is the number of contractions
of $T$ onto $S$. If $C_r$ denotes $r-$corona, then ${C_n \choose C_m} = {n \choose m}$.\\
Also\\
$$ \displaystyle \sum_S {T \choose S} ={ n \choose m}$$
if $deg(T) = n$ and the summation is extended over all $S$ with $deg(S) = m$ where the degree $deg(T)$ is the number of
leaves of $T$.\\
If $(1 + x)^T$ denotes the $T-th$ power of $ 1 + x$, then we obtain the nonassociative binomial formula\\
$(1 + x)^T = \sum {T \choose S} (x^S)$ where the summation is extended over all trees $S$.\\
The $T-th$ root $\sqrt[T]{1+x}$ is the unique power series $f$ with constant coefficient $1$ and $f^T = 1 + x$.\\
It is uniquely defined and its coefficients are denoted by $T^{-1} \choose S$ They are computed in section 6,
if $T$ is a corona.\\
It is an open problem to device a procedure for the computation of the coefficients $T^{-1} \choose S$
of $\sqrt[T]{1+x}.$\\
Also we consider binomial coefficients ${T \choose S, V}$ of the second kind which are defined to be the number of
contractions of $T$ onto $S$ for which the complementary contractions are onto $V$. It is shown that \\
$$(x \otimes 1 + a \otimes x)^T = \sum {T \choose S, V} x^S \otimes x^V$$

\medskip

\section{ Binomial coefficients of the first kind}

Denote by $\mathbb{P}$ the set of isomorphism classes of finite, planar, reduced rooted trees. The empty tree
$\underline{1}$ is included in $\mathbb{P}$. For  $ T \in \mathbb{P}$ and any subset $I$ of the set $L(T)$ of leaves of $T$ denote of
$T | I$ the contraction of $T$ onto    $I$,  see [G 3].\\
\begin{defn}
For any pair $T, S$ of $\mathbb{P}$ denote by $${T\choose S}$$
the number of $I \subseteq L(T)$ such that the contraction $T|I$ of $T$ onto $I$ is isomorphic with $S$.\\
We call $T \choose S$ the planar binomial coefficient of $T$ over $S$ of the first kind
\end{defn}
\begin{proposition}
The following elementary properties are easy to check:
\begin{itemize}
\item[(i)]
${T \choose S}  = 0,$ if $deg(S) > deg(T)$
\item [(ii)]
If $deg(S) = deg(T),$ then $ {T \choose S}  = 1$ if $S = T$ and ${T \choose S}  = 0$ of $S \not= T.$
\item [(iii)]
${T \choose 1} = 1$
\item [(iv)]
If $x$ denotes the tree with a single vertex, then
${T \choose x}  = deg(T)$
\item [(v)]
  ${T \choose x^2 } =  {n \choose 2}$, if $n=deg(T)$ and
if $x^2$ denotes the unique tree of degree $2$ in $\mathbb{P}$.
\item [(vi)]
$\displaystyle \sum_{S \in \mathbb{P},deg(S) = m} {T \choose S} =  deg {(T) \choose  m}$

\end{itemize}
\end{proposition}
Let $T_1,  ...
 , T_m \in \mathbb{P}$ and denote by $\cdot_m(T_1, ..., T_m)$ or $T_1 \cdot T_2  \cdots T_m$ the
$m$-ary grafting of $T_1, ... , T_m$. It is the unique tree $T$ in $\mathbb{P}$ for which $T- \rho_T$ is the planar forest
which is the disjoint union of the ordered systems $T_1 \dot{\bigcup}\ ...\ \dot{\bigcup}\   T_m $ where $T_i$ is less than $Tj$ in case $i$ is
 less than
$j$.Here $\rho_T$ is the root of $T$.\\

\begin{proposition}
Let $m, r \geq 2$ and $T=T_1 \cdot T_2 \cdot ... \cdot T_m$ and $S= S_1 \cdot ... \cdot S_r \in \mathbb{P}$ and $T_i
\not= \underline{1}, S_j \not= \underline{1}$ for all $i, j.$\\
Then we have the following recursion formula\\
$ {T \choose S} = \displaystyle \sum_{1 \leq i_1 < i_2 < .. < i_r \leq m} {T_{i_1} \choose S_1} \cdot ... \cdot$
${T_{i_r} \choose S_r} + \displaystyle \sum^m_{i=1}$
${T_i \choose S} $
\end{proposition}
\begin{example}
Let $T = T_1 \cdot ... \cdot T_m, deg (T_i) = n_i,$ then\\
${T \choose x x^2} = \displaystyle \sum^m_{i=1}  {Ti \choose xx^2} + \displaystyle \sum_{1 \leq i < j \leq m}$
 $ n_i$ $n_j \choose 2$\\
 ${T \choose x^2 x} = \displaystyle \sum^m_{i=1} {T_i \choose x^2 x} + \displaystyle \sum_{1 \leq i < j \leq m}$
 ${n_i \choose 2} n_j$\\
 and\\
 ${n \choose 3} = \displaystyle \sum^m_{i=1} {n_i \choose 3} + \displaystyle \sum_{1 \leq i < j \leq m} \frac{n_i\ n_j}{2}$
 $(n_i + n_j -2)$\\
 If $m = 2,$ then $n_1 + n_2 = n$\\ and\\
 ${n \choose 3} = {n_1 \choose 3} + {n_2 \choose 3} + \frac{n_1\ n_2} {2} (n_1 + n_2 - 2)$
\end{example}

 \medskip

 \section {Coefficients of the second kind}

 Let $T, S, V \in \mathbb{P}$.\\

 \begin{defn}
 The number \\
$T \choose S, V$ $: =\# \{ I \subseteq L (T)\colon T | I = S, T| (L(T) - I) = V\}$\\
is called the planar binomial coefficient of $T$ over $S$ and $V$ of the second kind.
\end{defn}
In combinatorics there are well-known Stirling numbers of the first and second kind which has been a reason for
chooring this terminology.\\
From the definition it immediately follows that
$\displaystyle \sum_{V \in \mathbb{P}} {T \choose S,} V  = { T \choose S }$\\
and that
${T \choose S, V} = { T \choose V, S} $

For $m, r, s \geq 2$ denote by $ \Gamma (m, r, s)$ the set of all pairs $ (\alpha, \beta)$ of subsets of\\
 $\underline{m} \colon = \{ 1, 2, ... , m \}$ with the following properties.\\
 \begin {itemize}
 \item[(i)]
$ \# \alpha \in \{1, r \}, \# \beta \in \{1, s\}$
\item [(ii)]
$\alpha \cup \beta = \{ 1, 2, ... , m\}$
 \end{itemize}
 Then $\Gamma^\prime ( m, r, s) \colon = \{(\alpha, \beta) \in \Gamma(m, r, s) \colon \# \alpha = 1\}$\\
 is equal to $\{(\{i\}, \underline{m}) \colon 1 \leq i \leq r\}$ if $s = m$ and is equal to
 $$\{(\{i\}, \underline{m} - \{i\}) \colon 1 \leq i \leq r\}$$
 if $s = m - 1$ and $\Gamma^\prime (m, r, s)$ is empty otherwise.\\
 Also $\Gamma''(m, r, s) \colon = \{(\alpha, \beta) \in \Gamma'' (m, r, s) \colon \# \beta = 1\}$\\
 is equal to
 $$ \{(\underline{m}, \{i\}) \colon 1 \leq i \leq m\}$$
 if $r = m$ and is equal to
 $$\{(\underline{m} - \{i\}) \colon 1 \leq i \leq m\}$$
 if $r = m - 1$ which it is empty otherwise.\\
 Let $\Gamma^\ast(m, r, s) \colon = \Gamma (m, r, s) - (\Gamma^\prime ( m, r, s) \cup
 \Gamma''(m, r, s)).$\\
Then
 $\Gamma^\ast (m, r, s)$ is empty, if $r + s < m$\\
and if
$r + s \geq m,$ then\\
$\# \Gamma^\ast(m, r, s) = $ $m \choose r$ $r \choose r + s-m$ $ = \frac{m\ !}{(m-r)\ !\ (r + s - m )\ !\  (m - s)\ !}$

\medskip

Let
 $T, S, V \in \mathbb{P}$ with $ar (T) = m, ar(S)= r$ \\
$ar(V) = s, T = T_1 \cdot ... \cdot T_m, S= S_1 \cdot ... \cdot S_r, V = V_1 \cdot ... \cdot V_s, T_i, S_i, V_i \in \mathbb{P}-\{\underline{1}\}$
and\\
 $\delta = ( \alpha, \beta) \in \Gamma (m, r, s)$\\
Then \\
${T \choose S, V}$
$_{\gamma}$ is defined to b the number of subsets $I$ of $L(T) = L (T_1)\  \dot{\bigcup} ...\  \dot{\bigcup}\
L (T_m)$\\
such that, if
$$\alpha \colon= \{i \in \underline{m} \colon I \cap L (T_i) \not= \phi\}$$
$$ \beta \colon= \{i \in \underline{m} \colon I^\prime \cap L (T_i) \not= \phi\}$$
where $I^\prime$ is the complement of $I$ in $L(T)$ Then
$\gamma = (\alpha, \beta) \in \Gamma (m, r, s)$ and $T\  | \ I = S,\ T\  |\  I^\prime = V.$\\

\begin{proposition}
$ T \choose S, V$  $=$ $\displaystyle \sum_{\gamma \in \Gamma (m, r, s)}$ $T \choose S, V$ $_{\gamma}$
\end{proposition}
\begin{proposition}
\begin{itemize}
\item[(i)]
Let $\delta \in \Gamma^\prime (m, r, s), \delta = (\{i\}, \underline{m} - \{i\})$. Then $s = m-1$
and\\
$T \choose S, V$ $_{\delta}$ $=$ $T_1 \choose \underline{1}, V_1$ $\cdot ... \cdot$ $ T_{i-1} \choose \underline{1}
, V_{i-1}.$ $\cdot$ $T_i \choose S, \underline{1}$ $\cdot$ $T_{i+1} \choose \underline{1}, V_i$ $\cdot ... \cdot$
$T_m \choose \underline{1}, V_{m-1}$
\item[(ii)]
Let $\delta \in \Gamma^\prime (m, r, s), \delta = (\{i\}, \underline{m}) \cdot$\\
Then $s = m$
 and\\
$T \choose S, V$ $_{\delta}$ $=$ $T_1 \choose \underline{1}, V_1$ $\cdot ... \cdot$ $T_i \choose S, V_i$
$T_{i+1} \choose \underline{1}, V_{i+1}$ $\cdot ... \cdot$ $T_m \choose \underline{1}, V_m$\\
\item [(iii)]
Let $\gamma \in \Gamma^\ast (m, r, s,),\gamma = (\alpha, \beta)$ Then $\# \alpha = r, \# \beta = s$ and
$T \choose S, V$
$_{\gamma}$ $= d_1 \cdot d_2 \cdot ... \cdot d_m$ \\
where $d_i = $ $T_i \choose S_j,  V_k$ if $i \in \alpha \cap \beta$\\
and $\alpha $
$(j) = i, \beta(k)= i$ there $\alpha (j)$ is the $j-th$ element of $\alpha$ and $\beta(k)$ is the $k-th$ element
of $\beta$ in the natural order on $\underline{m}.$\\
If $i \in \alpha, i \notin \beta,$ then
$d_i = $ $T_i \choose S_j, \underline{1}$
with $\alpha(j) = i$\\
If $i \in \beta, i \notin \alpha,$ then
$d_i=$ $T_i \choose \underline{1}, V_k$
with $\beta(k) = i.$
\end{itemize}
\end{proposition}
\section{Planar polynomials}
Let $K$ be a field and denote by $K\{x\}$ the $K-$algebra of planar polynomials.
For any $T \in \mathbb{P}, T \not= \underline{1},n= deg(T),$ there is a K- multilinear operation\\
$\cdot_T \colon (K\{x\})^n \rightarrow K \{x\}$ mapping $(S_1, ... , S_n)\in \mathbb{P}^n$ onto the T-grafting\\
$\cdot_T (S_1, ..., S_n)$ which is a planar tree $V$ in $\mathbb{P}$ containing $T$ such that $V - In(T)$ is the planar
forest $S_1 \dot{\bigcup} ... \dot{\bigcup} S_n$ there $In(T)$ denotes the set of inner vertices of $V$.
\begin{defn}
Let $f \in K \{x\}.$ Then
$$ f^T \colon = \cdot_T(f,..., f )$$
 is called the $T-th$ power of $F$
\end{defn}
\begin{proposition}
$(1 + x)^T = \displaystyle \sum_{S \in \mathbb{P}}$ ${T \choose S} x^S$
\end{proposition}
\begin{corollary}
For any $f \in K\{x\}$ we have $ \colon (1 +f)^T= \sum {T \choose S} f^S$
\end{corollary}
Denote by $\frac{d}{dx} \colon K\{x\} \rightarrow K \{x\}$ the derivation relative to $x$. It is a K-linear map satisfying
the usual product rule for $m-ary$ products such that $\frac{d}{dx}  (x) = 1$
\begin{proposition}
For any $T \in \mathbb{P},n = deg(T),\\
 \frac{d}{dx}(T(x)) = \displaystyle \sum_{deg(U)=n-1} {T \choose U} x^U$
\end{proposition}
\begin{proposition}
Let $T, S \in \mathbb{P}, T \not= 1 \not= S.$\\
Then for any $f \in K\{x\}$ we get \\
$$(f^T)^S = f^{(T^S)}$$
\end{proposition}
\begin{corollary}
$\biggl(\displaystyle \sum_{V \in \mathbb{P}}$ $T \choose V$ $x^V \biggl)^S= \sum$ $T^S \choose V$ $x^V$
\end{corollary}
\section{Co-addition}
Let $K\{x\} \otimes K \{x\}$ be the tensor product of the $K-$algebra $K\{x\}$ with itself. It is an algebra over
the operad $F$ freely generated by operation $(\mu_m)_{m \geq 2}$ where the degree of $\mu_m$ is $m$.\\
There is a unique $K-$ algebra homomorphism
$$\Delta \colon K \{x\} \rightarrow K\{x\} \otimes K\{x\}$$
such that $$\Delta(x) = x \otimes 1 + 1 \otimes x$$
\begin{proposition}
$\Delta (T(x)) = \displaystyle \sum_{S, V \in \mathbb{P}} {T \choose S,} x^S \otimes x^V$
\end{proposition}
\section{Planar roots}

Let $K\{\{x\}\}$ be the $K-$algebra of planar power series in $x$.
\begin{proposition}
Let $char (K) = 0$ and $T \in \mathbb{P}. $\\
There is a unique power series
$$ f \in K\{\{x\}\}$$
such that $f^T = (1+x)$ and whose constant term $f(0) = 1.$\\
We denote $f$ by $\sqrt[T]{1+x}$ or $(1 + x)^{T^{-1}}$
and call it the $T-th$ root of $1 + x$.
\end{proposition}

${\bf Example}$
Let $T =  x^2$ be the unique tree in $\mathbb{P}$ of degree 2. Then $\sqrt[T]{1+x}$ is denoted simply by $\sqrt{1+x}$.\\
One can show that
$\sqrt   {1+x}= \displaystyle \sum_{n=deg(S)>0\atop S\  \textrm{binary}} \frac{(-1)^{n-1}}{2^{2n-1}} x^S$\\
A tree $S$ in $\mathbb{P}$ is binary, if the arity $ar_S(v)$ of each inner vertex of $S$ is equal to $2$.\\
The arity $ar_S(v)$ is the number of edges of $S$ which are incident with $v$ and which are not on the path between $v$ and
the root of $S$.\\
\begin{remark}
One can show that
$$\sqrt[T]{1+x}= exp_T(\frac{1}{n} log_T(1 + x)$$
if $n= deg(T)$ and $exp_T$ is the exponential series relative to $T$ which means that
$$(exp_T(x))^T = exp_T(n x)$$
$$ ord(exp_T(x) - 1 - x) \geq 2$$
and $log_T$ is the inverse of $exp_T$, see [G1], [G2] for the definitions of $exp_T, log_T$.\\
\end{remark}
For any $T \in \mathbb{P} $, let $n_k(T)$ be the number of vertices of $T$ of arity $k$. Then
$n_0(T)$ is the degree of $T$ and $n_1(T)=0 $ is equivalent to $T$ being reduced.\\
If $T = T_1 \cdot ... \cdot T_m, T_i \in \mathbb{P}, T_i \not= \underline{1}$ for all $i$, then
$$n_k(T) = \displaystyle \sum^m_{i=1} n_k (T_i)$$
if $k \not= m$ and
$$ n_m(T) = 1+ \displaystyle \sum^m_{i=1} n_m (T_i).$$
One can show that
$$\displaystyle \sum_{k \geq 2} (k-1) n_k(T) = n_0(T)-1$$
Let $ m(T) = \displaystyle \sum^{< 00}_{k=0} n_k(T)$
\begin{proposition}
Let $C_d = x^d$ be the $d-corona, d \geq 2,$ and
$$\sqrt[C_d]{1 + x}$$
the $C_d-$ root of $(1 + x).$\\
Denote by ${C^{-1}_d \choose T}$ the coefficient of
$$\sqrt[C_d]{1+x}$$
relative to $T \in \mathbb{P}$.\\
Then $C_d^{-1} \choose T$ $ = 0$ if there is a vertex $v$ in $T$ with $ar_T(\vee) \geq d+1.$\\
If $ ar_T(v) \leq d$ for all vertices $v$ of $T$, then\\
${C^{-1}_d \choose T} = (-1)^{n(T) - n_0(T)} \cdot \displaystyle \prod_{k \geq 2} \biggl(\frac{1}{k} {d-1 \choose k - 1}\biggl)^{n_k(T)}$
\end{proposition}
Let $v= (\nu_k) k\geq 2$ be a sequence of integers $\nu_k \in \mathbb{N}$ such that $\nu_k=0$ for almost all $k$.\\
Let $deg(\nu)= \nu_0= 1 + \displaystyle \sum^\infty_{k=2} (k-1) k$ and $\mathbb{P}(\nu) \colon = \{T \in \mathbb{P}
\colon n_k (T) = \nu_k$ for all $k \geq 2\}$.\\
Then $\nu_0 = deg(T)$.\\
Denote by $C(\nu)$ the number of trees in $\mathbb{P}(\nu) j$ it is called Catalan number of $\nu$.\\
If $\nu_k=0$ for $k \geq 3,$ then $\mathbb{P}(\nu)$ is the set of binary trees of degree $\nu_0= 1 + \nu_2$ and $C (\nu)$
is the usual Catalan number.\\
Let $a_\nu \colon = (-1)^{\overline{\nu} - \nu_0} \displaystyle \prod_{k\geq 2}\biggl(\frac{1}{k} {\nu_{0 -1}
 \choose k-1}\biggl)^{\nu_k}$\\
 Let  ${C^{-1}_d \choose T} = a_{\nu (T)}$, if $\nu(T) = \biggl(n_k(T)\biggl)_{k\geq 2}$ by Proposition 6.3
 \begin{corollary}
 $\displaystyle \sum_{deg(\nu) = n} \# C(\nu) \cdot a_v = {1/d \choose n}= (-1)^{n-1} \frac{1}{n! d^n} \displaystyle
 \prod^{n-1}_{r=0} (rd-1)$
 \end{corollary}
 \begin{proof}
 The canonical algebra homomorphism
 $$\nu \colon K \{\{c\}\} \rightarrow K [[x]]$$
 mapping $T$ onto $x^{deg(T)}$ maps
 $(1+x)^{c^{-1}_d}$ onto $(1 + x)^{1/d} = \displaystyle \sum^\infty_{n=0}  {1/d \choose n} x^n$.
 \end{proof}

\textbf{Open problems:} If the coefficient of $\sqrt[T]{1+x}$ relative to $S \in \mathbb{P}$ is denote by $\overline{T}^1 \choose
S$, then the question arises by what procedure this binomial coefficient $\overline{T}^1 \choose S$ can be computed.\\
If there a recursive formula to compute $\overline{T}^1 \choose S$.\\
More generally there is a power series $f$ such that
$$f(0) = 1$$
$$f^T = (1 + x)^S$$
One could denote $f$ by $(1 + S)^{S\overline{T}^1}.$
How to compute the coefficients of $f$?

\end{document}